\documentclass[12pt, reqno]{amsart}
\usepackage{amsmath, amsthm, amscd, amsfonts, amssymb, graphicx, xcolor, cite}
\usepackage[bookmarksnumbered, colorlinks, plainpages]{hyperref}
\usepackage[foot]{amsaddr}
\newcommand{\pa}{\nabla}
\newtheorem{theorem}{Theorem}[section]

\newtheorem{corollary}[theorem]{Corollary}
\theoremstyle{definition}

\theoremstyle{remark}

\numberwithin{equation}{section}
\begin{document}
\setcounter{page}{1}

\centerline{}

\centerline{}

\title[]{Global Well-posedness for Incompressible Hookean Elastodynamics in the Critical Besov Spaces}
\author[ ]{Zexian Zhang$^{1}$ \and Yi Zhou$^{2}$}

\address{$^{1}$Shanghai Center for Mathematical Sciences, Fudan University, Shanghai 200433, China.}
\email{\textcolor[rgb]{0.00,0.00,0.84}{23110840019@m.fudan.edu.cn}}

\address{$^{2}$School of Mathematics Science, Fudan University, Shanghai 200433, China.}
\email{\textcolor[rgb]{0.00,0.00,0.84}{yizhou@fudan.edu.cn}}




\begin{abstract}
We identify the wave maps type nonlinearities of incompressible Hookean elastodynamics equations in Lagerangian coordinates, and iterate them in the adapted $U^2$-type spaces to prove the small data global well-posedness in the critical Besov space $\dot{B}^{\frac{n}{2}+1}_{2,1}(\mathbb{R}^n)\times \dot{B}^{\frac{n}{2}}_{2,1}(\mathbb{R}^n)\ (n\ge 2)$.
\end{abstract}

\maketitle

\section{Introduction}

The goal of this paper is to prove the global well-posedness for the Cauchy problem of incompressible isotropic Hookean elastodynamics with small initial data in  $\dot{B}^{\frac{n}{2}+1}_{2,1}(\mathbb{R}^n)\times \dot{B}^{\frac{n}{2}}_{2,1}(\mathbb{R}^n)\ (n\ge 2)$.

We start by fomulating the incompressible Hookean elastodynamics equations in Lagrangian coodinates. Let $y$ be the lagerangian coodinate, $X(t,y)$ be the flow map. We call $X(t,y)$ is incompressible, if for any connected domain $\Omega$, we have
$$\int_\Omega \mathrm{d}y = \int_{\Omega_t}\mathrm{d}X,\ \Omega_t = \{X(t,y)|\ y\in \Omega\}.$$

Clearly the incompressibilty is equivalent to \begin{equation}\label{con1}
  \det(\nabla X) \equiv 1.
\end{equation}

In elastodynamics, the motion of the elastic
fluids for the homogeneous, isotropic, and hyperelastic materials in Lagrangian coordinates is determined by the following Lagrangian functional
of flow maps:
\begin{equation} \label{la}
    \begin{aligned}
        \mathcal{L}(X,T)=&\int_{\mathbb{R}^{1+n}}
        \biggl( \frac{1}{2}| \partial_tX(t,y)|^2 - W(\nabla X(t,y))\\
        &+ p(t,y)\bigl[\det\bigl(\nabla X(t,y)\bigr)- 1\bigr] \biggr)
        \mathrm{d}y\mathrm{d}t,
    \end{aligned}
\end{equation}
where $W \in C^\infty (\text{GL}_n(\mathbb{R}),\mathbb{R}^+)$ is the strain energy function 
depending on the deformation tensor $\nabla X$, and pressure $p(t,y)$ is a Lagrangian
multiplier which ensures the flow maps to be incompressible.

In the Hookean case, the strain energy function is given by 
\begin{equation}
  W(\nabla X) = \frac{1}{2}|\nabla X|^2.
\end{equation}

To proceed further, denote
$$X(t,y) = Y(t,y)+y.$$

Then the Euler-Lagrangian equation of \eqref{la} takes the form of 
\begin{equation}\label{ela1}
  \Box Y = -\nabla p - (\nabla Y)^{T}\Box Y.
\end{equation}
Here and in what follows, we write $\Box  = \partial_t^2-\Delta$ for the D'Alembertian, and $(\nabla Y)_{ij} = \nabla_j Y_i$ for the Jacobian of $Y$.

The constraint \eqref{con1} becomes
\begin{equation}\label{ela2}
  \nabla \cdot Y = -\sum_{k=2}^n E_k(\nabla Y).
\end{equation}
Here $E_k(A)$ denotes the sum of principal minors of size $k$ in the $n\times n$ matrix $A$, and we use the algebraic identity
$$\det(tI-A) = t^n+\sum_{k=0}^{n-1}(-1)^{n-k}E_{n-k}(A)t^k.$$

Next we show that \eqref{ela1} and \eqref{ela2} are actually a hyperbolic-elliptic coupled system. Let $\bf{P} = 1- \Delta^{-1}(\nabla \text{div})$ be the Leray projector. We split
\begin{equation*}
  Y = Y^{\text{df}}+Y^{\text{cf}},
\end{equation*}
where $Y^{\text{df}} = \mathbf{P}Y$ is the divergence-free part and $Y^{\text{cf}} = (1-\mathbf{P})Y$ is the curl-free part.

Applying $\mathbf{P}$ to \eqref{ela1}, we obtain 
\begin{equation}\label{ela3}
  \begin{aligned}
  \Box Y^{\text{df}}_j &= -\mathbf{P}_j[(\nabla Y)^T\Box Y]  \\
  & = -\Delta^{-1}\pa_k [\pa_k ((\nabla Y)_{lj}\Box Y_l) - \pa_j ((\nabla Y)_{lk}\Box Y_l)]\\
  & = -\Delta^{-1}\pa_k [\pa_k (\nabla_j Y_l\Box Y_l) - \pa_j (\nabla_k Y_l\Box Y_l)]\\
  & = -\Delta^{-1}\pa_k [(\nabla_j Y_l)\Box \pa_k Y_l - (\nabla_k Y_l)\Box \pa_j Y_l].
  \end{aligned}
\end{equation}
Here we use the summation convention over repeated indices.

And the constraint \eqref{ela2} turns into
\begin{equation}\label{ela4}
  Y^{\text{cf}}_j = \Delta^{-1}\nabla_j(\nabla \cdot Y^{\text{cf}}) = -\Delta^{-1}\nabla_j \sum_{k=2}^n E_k(\nabla Y).
\end{equation}

Combining above, we refomulate the incompressible Hookean elastodynamics equations into the following hyperbolic-elliptic coupled system
\begin{equation}\label{ela5}
    \begin{cases} \Box \nabla_i Y_j^{\text{df}} = R_i R_k [(\nabla_j Y_l)\Box \pa_k Y_l - (\nabla_k Y_l)\Box \pa_j Y_l] ,  \\ \nabla_i Y^{\text{cf}}_j = R_i R_j \sum_{k=2}^n E_k(\nabla Y), \\
    (Y, \partial_t Y)|_{t=0} =(f,g),  \end{cases}
  \end{equation} 
where $R_i = (-\Delta)^{-\frac{1}{2}}\pa _i$ is the Riesz-type operator, and the initial data admits the compatibility condition
$$\det(I+\nabla f) =1,\ \text{Tr}((I+\nabla f)^{-T}\nabla g) = 0.$$

It's clear that the equation \eqref{ela5} is invariant under the scaling 
$$Y(t,y)\mapsto \lambda^{-1}Y(\lambda t,\lambda y),$$
and standard scaling analysis applied to the norm  $\left\| Y \right\|_{\dot{H}^s(\mathbb{R}^n)}$ shows that the critical Sobolev exponent is $s=\frac{n}{2}+1$.

For the global existence of the incompressible elastodynamics with the small data of high regularity, the three-dimensional case is proved by Sideris-Thomas  in \cite{sideris_Global_2005} and \cite{sideris_Global_2007} in Eulerian coordinates, and in two-dimensional case the almost global existence is established by Lei-Sideris-Zhou \cite{lei_Almost_2015} in Eulerian coordinates. Later this is improved to global existence by Lei \cite{lei_Global_2016} in Lagrangian coordiantes. The same results are established independently by Wang \cite{wang_Global_2017}  and Cai \cite{cai_Vanishing_2019} in Eulerian coordinates.  For more related topics and results we refer the reader to \cite{cai_Global_2022}, \cite{cai_Uniform_2023}, \cite{lei_Global_2008}, \cite{lei_Global_2005}, \cite{lin_Hydrodynamics_2005} and \cite{zhu_Global_2018}. 

In this paper we focus on the small data result at the critical regularity. More precisely, we will prove
\begin{theorem}\label{main}
  There exists $\varepsilon >0$ such that for any initial data satisfying $$\left\| (f,g) \right\|_{\dot{B}^{\frac{n}{2}+1}_{2,1}\times \dot{B}^{\frac{n}{2}}_{2,1}}\le \varepsilon,$$
  there exists a unique global solution to \eqref{ela5} satisfying
  $$\left\| (\nabla Y, \partial_t \nabla Y) \right\|_{L^\infty_t (\dot{B}^{\frac{n}{2}}_{2,1}\times \dot{B}^{\frac{n}{2}-1}_{2,1})}\lesssim  \varepsilon.$$
Moreover, the solution depends continously on the initial data. 
\end{theorem}
Considering the nonlinearities of the \eqref{ela5}, it's sufficient to find a scaling invariant iteration space $S^{\frac{n}{2}}$ with the properties
\begin{equation}\label{pro}
  S^{\frac{n}{2}}\cdot S^{\frac{n}{2}} \hookrightarrow S^{\frac{n}{2}}\ \text{and}\ \Box^{-1}(S^{\frac{n}{2}}\Box S^{\frac{n}{2}})\hookrightarrow S^{\frac{n}{2}},
\end{equation}
see Theorem \ref{S1} below for the precise statement.

Such a space has been constructed to tackle the wave maps equations. To solve the problem in the critical Besov space $\dot{B}^{\frac{n}{2}}_{2,1}\times \dot{B}^{\frac{n}{2}-1}_{2,1}$, Tataru \cite{tataru_Global_2001}  introduces null frames and characteristic energy estimates, and Candy-Herr \cite{candy_Division_2018}  introduces the adapted $U^2$ type atomic space and $L^2_{t,x}$ bilinear estimates. The present paper will borrow quite heavily form  \cite{candy_Division_2018}, and in fact, will be built on the   adpated $U^2$ function spaces and bilinear estimates established there.

In the critical Sobolev space, the properties \eqref{pro} may not be appropriate because of logarithmic divergences arising from certain types of low-high interactions. For the wave maps equation, the problem is solved via renormalization argument. The sphere target case was solved by Tao \cite{tao_Global_2001},
later by Krieger \cite{krieger_Global_2004} for the hyperbolic plane target, and for more general targets by
Klainerman-Rodnianski \cite{klainerman_Global_2001} and Tataru \cite{tataru_Rough_2005}.

The remaining
part of this paper is organized as follows: In Section 2 we will introduce the basic notations. In Section 3 we will follow \cite{candy_Division_2018} to introduce the iteration space based on the adapted $U^{2}$ space and summarize related properties and estiamtes. In Section 4 we will give the proof of Theorem \ref{main} via the estiamtes shown in Section 3.
\section{Notations}
We use $A\lesssim B $ to denote the statement that $A \le  CB$ for some absolute constant $C$, and $A \sim  B$ to denote the statement $A\lesssim B\lesssim A$.

Given funtions $f(x)$ on $\mathbb{R}^{n}$ and $u(t,x)$ on $\mathbb{R}^{1+n}$, we use $\hat{f}(\xi)$ to denote its spatial Fourier tranformation, and $\tilde{u}(\tau,\xi)$ to denote its space-time Fourier tranformation. Let $|\nabla|$ be the Fourier multiplier such that $\widehat{|\nabla|f}(\xi) = |\xi| \hat{f}(\xi)$.

For dyadic numbers $\lambda \in  2^{\mathbb{Z}}$, let $P_\lambda$ denote a smooth spatial cutoff to the Fourier region $|\xi|\approx \lambda$. We shall write $u_{\lambda}$ for $P_\lambda u$.

We use 
\begin{equation}
  V(t)(f,g) = \cos(t|\nabla|)f + \sin(t|\nabla|)|\nabla|^{-1}g 
\end{equation} to denote the solution to homogeneous wave eqaution, and
\begin{equation}
  \Box^{-1}F(t) = \mathbf{1}_{\{t\ge 0\}}\int^t_0 \sin((t-s)|\nabla|)|\nabla|^{-1}F(s) \mathrm{d}s
\end{equation} 
to denote the solution for inhomogeneous wave equation with vanishing initial data.

\section{Function spaces}
In this section we follow \cite{candy_Division_2018} to introduce the iteration space $S^{\frac{n}{2}}$ and corresponding estmates in preperation for the proof of theorem \ref{main}.  

\subsection{The $U^p$ and $V^p$ spaces}
For a  partition $\tau = (t_i)_{i=1}^N$, let 
$\mathcal{I}_{\tau} = \{[t_i,t_{i+1})\}_{i=1}^N$, where $t_{N+1} = \infty $. Let $\mathbf{T}$ be the set of all the partitions defined above. Let $1\le p<+\infty$. We say $u$ is a $U^p$-atom if there exists $\tau \in \mathbf{T}$ such that 
\begin{equation}
  u = \sum_{I \in \mathcal{I}_{\tau}}f_{I}\mathbf{1}_I(t),\ \ \bigl(\sum_{I \in \mathcal{I}_{\tau}}\left\| f_I \right\|^p_{L_{2}}\bigr)^{\frac{1}{p}}=1.
\end{equation}

We define the $U^p$ space to be the atomic space spanned by the $U^p$-atoms, with the norm
\begin{equation}
  \left\| u \right\|_{U^p} = \inf\{\sum_{i=1}^{\infty}|c_i|;\ u = \sum_{i=1}^\infty c_i u_i,\ u_i\ \text{are}\ U^p\ \text{atoms} \}.
\end{equation}

From the definition, functions $u$ in $U^p$ are bounded, have one-sided limits everywhere and right-continous with $\lim_{t \to -\infty}u=0$ in $L^2(\mathbb{R}^n)$.

The closely related $V^p$ space is defined to be the union of all right-continous function $v$ such that $ \lim_{t \to -\infty}v(t) = 0$ in $L^2(\mathbb{R}^n)$, with norm
\begin{equation}
  \left\| v \right\|_{V^p}  = \sup_{\tau \in \mathbf{T}}\bigg(\sum_{i=1}^{N-1}\left\| v(t_{j+1})-v(t_{j}) \right\|_{L^2}^p\bigg)^{\frac{1}{p}}.
\end{equation}

\subsection{The iteration space and related properties} The solution space is based on the adapted version of $U^2$ space. We define $U_{\pm}^2$ with the norm 
\begin{equation}
  \left\| u \right\|_{U_{\pm}^2} = \| e^{\pm it|\nabla|}u \|_{U^2}.
\end{equation}

Let the space $S$ be the subspace for functions $u\in C_t L^2_x$ with $|\nabla|^{-1}\partial_t u \in C_t L^2_x$, with the norm  
\begin{equation}
  \left\| u \right\|_S = \left\| u+i|\nabla|^{-1}\partial_t u \right\|_{U^2_+}+ \left\| u-i|\nabla|^{-1}\partial_t u \right\|_{U^2_-}.
\end{equation}

The solution sapce $S^{\frac{n}{2}}$ is defined as the subspace of $C_t \dot{B}^{\frac{n}{2}}_{2,1}$ with the norm
\begin{equation}
  \left\| u \right\|_{S^{\frac{n}{2}}} = \sum_{\lambda \in 2^{\mathbb{Z}}}\lambda^{\frac{n}{2}}\left\| u_\lambda \right\|_S,
\end{equation}
and it is automatically banach space since the subspace of continous functions in $U^{2}_{\pm}$ is complete. Moreover, it's clear that the Riesz-type operator $R_i$ is bounded in the iteration space $S^{\frac{n}{2}}$, i.e.
\begin{equation}\label{B1}
  \left\| R_i u \right\|_{S^{\frac{n}{2}}}\lesssim \left\| u \right\|_{S^{\frac{n}{2}}} .
\end{equation} 

Noticing the identities
\begin{equation}
  \begin{aligned}
    u &= \frac{1}{2}(u+i|\nabla|^{-1}\partial_t u)+\frac{1}{2}(u-i|\nabla|^{-1}\partial_t u), \\
    |\nabla|^{-1}\partial_t u& = \frac{1}{2i}(u+i|\nabla|^{-1}\partial_t u)-\frac{1}{2i}(u-i|\nabla|^{-1}\partial_t u),
\end{aligned}
\end{equation}

then we have
\begin{equation}
  \left\| (u_\lambda,\partial_t u_\lambda) \right\|_{L^\infty_t(L^2_x \times \dot{H}^{-1}_x)}\lesssim \left\| u_\lambda \right\|_S,\ \ \left\| (u,\partial_t u) \right\|_{L^\infty_t (B^{\frac{n}{2}}_{2,1}\times B^{\frac{n}{2}-1}_{2,1})}\lesssim \left\| u \right\|_{S^{\frac{n}{2}}}.
\end{equation}

To describe the space $S$ via the dual pairing $(DU^2)^* = V^2$, we need a weaker version of the space $S$ based on $V_{\pm}^2$. Let $S_w$ be the space consisting right-continuous functions with norm
\begin{equation}
  \left\| v \right\|_{S_w} = \left\| v \right\|_{V^2_+ + V^2_{-}} = \inf_{v = v_- +v_+} \left\| v_-  \right\|_{V_-} +\left\| v_+ \right\|_{V_+},
\end{equation}
where $\left\| v \right\|_{V_{\pm}^2} = \left\| e^{\pm it|\nabla|}v \right\|_{V^2}.$ Since $V^2 \hookrightarrow U^2$, we obtain
$$ \left\| u \right\|_{S_w}+ \left\||\nabla|^{-1}\partial_t u \right\|_{S_w}\le \left\| u+i|\nabla|^{-1}\partial_t u \right\|_{V^2_+}+ \left\| u-i|\nabla|^{-1}\partial_t u \right\|_{V^2_-}\lesssim \left\| u \right\|_S.$$

We summarize the main properties of the space $S$ in the following Theorem, and its proof can be found in \cite[Lemma 2.1, Lemma 3.1 and Theorem 3.2]{candy_Division_2018}.

\begin{theorem}\label{S1}
    \begin{enumerate}
        \item {\rm (Characteristic of $S$)} Let $(\psi,\partial_t \psi) \in C_b(\mathbb{R},L^2_x \times \dot{H}^{-1}_x)$ with $\lim_{t \to -\infty}\left\| (\psi,\partial_t \psi) \right\|_{L^\infty_t(L^2_x \times \dot{H}^1_x) }=0$, then
        \begin{equation}
          \begin{aligned}\left\| \psi \right\|_S\lesssim \sup _{{\tiny \begin{array}{c} \phi \in C^\infty _0 \\ \| \phi \| _{S_w} \leqslant 1 \end{array}}} \Big | \int _{\mathbb {R}} \langle \Box \phi , |\nabla |^{-1} \psi \rangle _{L^2_x} \mathrm{d}t \Big |\end{aligned},
        \end{equation}
whenerver the right hand side is finite.
         \item  Let $\chi \in C^\infty$ such that $\chi(t)\equiv 1$ for $t\ge 0$ and $\chi(t)\equiv 0$ for $t\le -1$. For $f_\lambda,g_\lambda \in L^2$, we have \begin{equation}
           \left\| \chi(t|\nabla|)V(t)(f_\lambda,g_\lambda) \right\|_{S}\lesssim \left\| f_\lambda \right\|_{L^2}+ \lambda^{-1}\left\| g_\lambda \right\|_{L^2}.
         \end{equation}
         \item  Given dyadic numbers $\lambda_0,\lambda_1,\lambda_2$. If $u_{\lambda_1},v_{\lambda_2}\in S$, then $P_{\lambda_0}(u_{\lambda_1}v_{\lambda_2}) \in S$ and 
         \begin{equation}
           \lambda_0^{\frac{n}{2}}\left\| P_{\lambda_0}(u_{\lambda_1}v_{\lambda_2}) \right\|_{S} \lesssim(\lambda_1 \lambda_2)^{\frac{n}{2}}\left\| u \right\|_{S}  \left\| v \right\|_{S} .
         \end{equation}
         \item Given dyadic numbers $\lambda_0,\lambda_1,\lambda_2 $. If $u_{\lambda_1},v_{\lambda_2} \in S $ and $\Box v_\lambda \in L^1_{t,\rm{loc}}L^2_x$, then $\Box^{-1}P_{\lambda_0}(u_{\lambda_1}\Box v_{\lambda_2}) \in S$ and we have 
         \begin{equation}
           \lambda_0^{\frac{n}{2}} \left\| \Box^{-1}P_{\lambda_0}(u_{\lambda_1}\Box v_{\lambda_2})\right\|_S \lesssim(\lambda_1\lambda_2)^{\frac{n}{2}}\left\| u_{\lambda_1} \right\|_S \left\| v_{\lambda_2} \right\|_S.
         \end{equation}
    \end{enumerate}
\end{theorem}

After summation over the frequency pieces in Theorem \ref{S1}, we have the following collary 
\begin{corollary}\label{S2} 
  \begin{enumerate}
    \item {\rm (Control of free waves)}  Let $\chi \in C^\infty$ such that $\chi(t)\equiv 1$ for $t\ge 0$ and $\chi(t)\equiv 0$ for $t\le -1$. For $(f,g) \in \dot{B}^{\frac{n}{2}}_{2,1}\times \dot{B}^{\frac{n}{2}-1}_{2,1}$, we have 
    \begin{equation}\label{B4}
      \left\| \chi(t|\nabla|)V(t)(f,g) \right\|_{S^{\frac{n}{2}}}\lesssim \left\| (f,g) \right\|_{\dot{B}^{\frac{n}{2}}_{2,1}\times \dot{B}^{\frac{n}{2}-1}_{2,1}}.
    \end{equation}
    \item {\rm(Algebraic property)} Given $u,v \in S^{\frac{n}{2}}$, we have
    \begin{equation}\label{B2}
      \left\| uv \right\|_{S^{\frac{n}{2}}}\lesssim \left\| u \right\|_{S^{\frac{n}{2}}}\left\| v \right\|_{S^{\frac{n}{2}}}.
    \end{equation}
    \item{\rm(Null form estiamte)}  Given $u,v \in S^{\frac{n}{2}}$, if $\Box v_{\lambda} \in L^1_{t,\text{loc}}L^2_{x}$ for any dyadic number $\lambda$, then $\Box^{-1}(u\Box v) \in S^{\frac{n}{2}}$ and we have 
    \begin{equation}\label{B3}
      \left\| \Box^{-1}(u\Box v) \right\|_{S^{\frac{n}{2}}}\lesssim \left\| u \right\|_{S^{\frac{n}{2}}}\left\| v \right\|_{S^{\frac{n}{2}}}.
    \end{equation}   
        \end{enumerate}
\end{corollary}
\section{Proof of Theorem \ref{main}}
\begin{proof}
  We use the standard fix point argument to construct our solution in $S^{\frac{n}{2}}$. For the working space, set
  $$S_0^\frac{n}{2} = \{u \in S^\frac{n}{2}| \forall \lambda \in 2^{\mathbb{Z}},\ \Box u_\lambda \in L^1_{t,\text{loc}}L^2_x\}$$
to ensure the nolinear expressions are well defined in $S$. Define the map $T: S_0^{\frac{n}{2}}\to S_0^{\frac{n}{2}}$ to be 
\begin{equation}
  \begin{aligned}
    T[\nabla Y ]_{ij}(t) &= \chi(t|\nabla|)V(t)(\nabla_i \mathbf{P}_j f,\nabla_i \mathbf{P}_j g)- \Box^{-1}\nabla_i \mathbf{P}_j[(\nabla Y)^T\Box Y] \\
    &\quad  +R_i R_j\sum_{k=2}^n E_k(\nabla Y).
  \end{aligned}
\end{equation}
It's clear that the fix point of $T$ is a solution on $[0,\infty)$. 

Collary \ref{S2} shows that
$$
\begin{aligned}
  \left\| \Box^{-1}\nabla_i \mathbf{P}_j[(\nabla Y)^T\Box Y]  \right\|_{S^{\frac{n}{2}}}&= \left\| \Box^{-1}R_iR_k [(\nabla_j Y_l)\Box \pa_k Y_l - (\nabla_k Y_l)\Box \pa_j Y_l] \right\|_{S^{\frac{n}{2}}}\\
  &\lesssim \left\| \Box^{-1} [(\nabla_j Y_l)\Box \pa_k Y_l - (\nabla_k Y_l)\Box \pa_j Y_l] \right\|_{S^{\frac{n}{2}}},\\
  &\lesssim \left\| \nabla Y\right\|_{S^{\frac{n}{2}}}^2.\\
  \| R_i R_j \sum_{k=2}^n E_k(\nabla Y)\|_{S^{\frac{n}{2}}}&\lesssim (1+\left\|\nabla Y \right\|^{n-2}_{S^{\frac{n}{2}}})\left\|\nabla Y \right\|^{2}_{S^{\frac{n}{2}}}.
\end{aligned}$$
Here we use the null form estimate \eqref{B3} and the boundedness of Riesz-type operator \eqref{B1} in the fisrt inequality, and the algebraic property \eqref{B2} in the second inequality. Together with the control of free waves \eqref{B4}, we obtain
$$
\begin{aligned}
\left\| T[\nabla Y] \right\|_{S^{\frac{n}{2}}}&\lesssim \left\| (f,g) \right\|_{\dot{B}^{\frac{n}{2}+1}_{2,1}\times \dot{B}^{\frac{n}{2}}_{2,1}} + (1+\left\|\nabla Y \right\|^{n-2}_{S^{\frac{n}{2}}})\left\|\nabla Y \right\|^{2}_{S^{\frac{n}{2}}}, \\
\left\| T[\nabla Y - \nabla Y'] \right\|_{S^{\frac{n}{2}}}&\lesssim (\left\|\nabla Y\right\|_{S^{\frac{n}{2}}}+\left\|\nabla Y\right\|_{S^{\frac{n}{2}}}^{n-1}+\left\|\nabla Y'\right\|_{S^{\frac{n}{2}}}+\left\|\nabla Y'\right\|_{S^{\frac{n}{2}}}^{n-1} )\\
&\quad\ \cdot \left\| \nabla Y - \nabla Y' \right\|_{S^{\frac{n}{2}}}.
\end{aligned}$$
If we choose $\varepsilon\ll1$, then there exists $C_1>0$ such that $T$ is a contraction in $\{\nabla Y \in S^\frac{n}{2}_0|\ \|\nabla Y\|_{S^{\frac{n}{2}}}\le C_1 \varepsilon\}$. This implies the existence and uniqueness of a fix point in banach space $S^{\frac{n}{2}}$. 

Given two solutions $Y^1,Y^2$ with initial data $(f_1,g_1)$ and $(f_2,g_2)$ respectively, we have
\begin{align*}
  (\nabla_i Y_j^{1} - \nabla_i Y_j^{2})(t) &= \chi(t|\nabla|)V(t)(\nabla_i\mathbf{P}_j(f_1-f_2),\nabla_i\mathbf{P}_j(g_1-g_2))\\
  &\quad -\Box^{-1}\nabla_i \mathbf{P}_j[(\nabla Y^1-\nabla Y^2)^T\Box Y^1]\\
  &\quad  -\Box^{-1}\nabla_i \mathbf{P}_j[(\nabla Y^2)^T\Box (Y^1-Y^2)]\\
  &\quad  +R_i R_j\sum_{k=2}^n E_k(\nabla Y^1)-E_k(\nabla Y^2).
\end{align*}

Following the similar procedure we have 
\begin{align*}
  \left\|\nabla Y^{1} - \nabla Y^{2}\right\|_{S^{\frac{n}{2}}}&\lesssim \left\| (f_1-f_2,g_1-g_2) \right\|_{\dot{B}^{\frac{n}{2}+1}_{2,1}\times \dot{B}^{\frac{n}{2}}_{2,1}}\\
  &\quad +  (\left\|\nabla Y^1\right\|_{S^{\frac{n}{2}}}+\left\|\nabla Y^1\right\|_{S^{\frac{n}{2}}}^{n-1}+\left\|\nabla Y^2\right\|_{S^{\frac{n}{2}}}+\left\|\nabla Y^2\right\|_{S^{\frac{n}{2}}}^{n-1} )\\
  &\quad\quad \cdot \left\| \nabla Y^1- \nabla Y^2 \right\|_{S^{\frac{n}{2}}}\\
  &\lesssim \left\| (f_1-f_2,g_1-g_2) \right\|_{\dot{B}^{\frac{n}{2}+1}_{2,1}\times \dot{B}^{\frac{n}{2}}_{2,1}} + \varepsilon \left\| \nabla Y^1- \nabla Y^2 \right\|_{S^{\frac{n}{2}}}.
\end{align*}

Since $\varepsilon\ll 1$, it's clear that
$$ \left\|\nabla Y^{1} - \nabla Y^{2}\right\|_{S^{\frac{n}{2}}}\lesssim \left\| (f_1-f_2,g_1-g_2) \right\|_{\dot{B}^{\frac{n}{2}+1}_{2,1}\times \dot{B}^{\frac{n}{2}}_{2,1}}.$$
This gives the continous denpendence on the initial data.
\end{proof}

{\bf Acknowledgement.} This work was supported by the National Natural Science Foundation of China [No. 12171097], the
Key Laboratory of Mathematics for Nonlinear Sciences (Fudan University), the Ministry of Education
of China, Shanghai Key Laboratory for Contemporary Applied Mathematics and Shanghai Science and
Technology Program [No. 21JC1400600].
\bibliography{ref}

\begin{thebibliography}{10}

\bibitem{cai_Global_2022}
Yuan Cai.
\newblock Global vanishing viscosity limit for the two dimensional incompressible viscoelasticity in lagrangian coordinates.
\newblock {\em Calc. Var.}, 61(3):93, June 2022.

\bibitem{cai_Uniform_2023}
Yuan Cai.
\newblock Uniform bound of the highest-order energy of the 2d incompressible elastodynamics.
\newblock {\em SIAM J. Math. Anal.}, 55(5):5893--5918, October 2023.

\bibitem{cai_Vanishing_2019}
Yuan Cai, Zhen Lei, Fanghua Lin, and Nader Masmoudi.
\newblock Vanishing viscosity limit for incompressible viscoelasticity in two dimensions.
\newblock {\em Comm Pure Appl Math}, 72(10):2063--2120, October 2019.

\bibitem{candy_Division_2018}
Timothy Candy and Sebastian Herr.
\newblock On the division problem for the wave maps equation.
\newblock {\em Ann. PDE}, 4(2):17, December 2018.

\bibitem{klainerman_Global_2001}
Sergiu Klainerman and Igor Rodnianski.
\newblock On the global regularity of wave maps in the critical sobolev norm.
\newblock {\em International Mathematics Research Notices}, 2001(13):655--677, January 2001.

\bibitem{krieger_Global_2004}
Joachim Krieger.
\newblock Global regularity of wave maps from \$r{\textasciicircum}(2+1)\$ to \$h{\textasciicircum}2\$. small energy.
\newblock {\em Commun. Math. Phys.}, 250(3):507--580, October 2004.

\bibitem{lei_Global_2016}
Zhen Lei.
\newblock Global well-posedness of incompressible elastodynamics in two dimensions.
\newblock {\em Comm Pure Appl Math}, 69(11):2072--2106, November 2016.

\bibitem{lei_Global_2008}
Zhen Lei, Chun Liu, and Yi~Zhou.
\newblock Global solutions for incompressible viscoelastic fluids.
\newblock {\em Arch Rational Mech Anal}, 188(3):371--398, June 2008.

\bibitem{lei_Almost_2015}
Zhen Lei, Thomas Sideris, and Yi~Zhou.
\newblock Almost global existence for 2-d incompressible isotropic elastodynamics.
\newblock {\em Trans. Amer. Math. Soc.}, 367(11):8175--8197, April 2015.

\bibitem{lei_Global_2005}
Zhen Lei and Yi~Zhou.
\newblock Global existence of classical solutions for the two-dimensional oldroyd model via the incompressible limit.
\newblock {\em SIAM J. Math. Anal.}, 37(3):797--814, January 2005.

\bibitem{lin_Hydrodynamics_2005}
Fang-Hua Lin, Chun Liu, and Ping Zhang.
\newblock On hydrodynamics of viscoelastic fluids.
\newblock {\em Comm Pure Appl Math}, 58(11):1437--1471, November 2005.

\bibitem{sideris_Global_2005}
Thomas~C. Sideris and Becca Thomases.
\newblock Global existence for three-dimensional incompressible isotropic elastodynamics via the incompressible limit.
\newblock {\em Comm Pure Appl Math}, 58(6):750--788, June 2005.

\bibitem{sideris_Global_2007}
Thomas~C. Sideris and Becca Thomases.
\newblock Global existence for three-dimensional incompressible isotropic elastodynamics.
\newblock {\em Comm Pure Appl Math}, 60(12):1707--1730, December 2007.

\bibitem{tao_Global_2001}
Terence Tao.
\newblock Global regularity of wave maps ii. small energy in two dimensions:.
\newblock {\em Commun. Math. Phys.}, 224(2):443--544, December 2001.

\bibitem{tataru_Global_2001}
Daniel Tataru.
\newblock On global existence and scattering for the wave maps equation.
\newblock {\em American Journal of Mathematics}, 123(1):37--77, 2001.

\bibitem{tataru_Rough_2005}
Daniel Tataru.
\newblock Rough solutions for the wave maps equation.
\newblock {\em American Journal of Mathematics}, 127(2):293--377, 2005.

\bibitem{wang_Global_2017}
Xuecheng Wang.
\newblock Global existence for the 2d incompressible isotropic elastodynamics for small initial data.
\newblock {\em Ann. Henri Poincar{\'e}}, 18(4):1213--1267, April 2017.

\bibitem{zhu_Global_2018}
Yi~Zhu.
\newblock Global small solutions of 3d incompressible oldroyd-b model without damping mechanism.
\newblock {\em Journal of Functional Analysis}, 274(7):2039--2060, April 2018.

\end{thebibliography}
\bibliographystyle{plain}

\end{document}